\documentclass[12pt]{amsart}
\usepackage{amssymb}
\newtheorem{theorem}{Theorem}
\newtheorem{proposition}[theorem]{Proposition}

\newtheorem{lemma}[theorem]{Lemma}
\newtheorem{example}[theorem]{Example}

\title{Remarks on Lempert functions of balanced domains}

\author{Nikolai Nikolov and Peter Pflug}

\address
{Institute of Mathematics and Informatics\\ Bulgarian Academy of
Sciences\\ Acad. G. Bonchev 8, 1113 Sofia,
Bulgaria}\email{nik@math.bas.bg}

\address{Carl von Ossietzky Universit\"at Oldenburg\\
Institut f\"ur Mathematik, Fakult\"at V\\ Postfach 2503\\ D-26111
Oldenburg, Germany}\email{pflug@mathematik.uni-oldenburg.de}



\begin{document}

\begin{thanks}{This note was written during the stay of the first
named author at the Universit\"at Oldenburg supported by a grant
from the DFG, Az.~PF 227/9-1 (October 2007). He thanks both
institutions for their support.}
\end{thanks}

\begin{abstract}
This note should clarify how the behavior of certain invariant
objects reflects the geometric convexity of balanced domains.
 \end {abstract}

\maketitle

\section{Introduction and results}

By $\Bbb D$ we denote the unit disc in $\Bbb C.$ Let $D$ be a
domain in $\Bbb C^n.$ Recall first the definitions of the
Carath\'eodory pseudodistance and the Lempert function of $D:$
$$c_D(z,w)=\sup\{\tanh^{-1}|f(w)|:f\in\mathcal O(D,
\Bbb D):f(z)=0\},$$
$$\tilde
k_D(z,w)=\inf\{\tanh^{-1}|\alpha|:\exists\varphi\in\mathcal O(\Bbb
D,D):\varphi(0)=z,\varphi(\alpha)=w\}.$$

The Kobayashi pseudodistance $k_D$ can be defined as the largest
pseudodistance below of $\tilde k_D$. Note that if $k_D^{(m)}$
denotes the $m$-th Lempert function of $D$, $m\in\Bbb N$, that is,
$$k_D^{(m)}(z,w)=\inf\{\sum_{j=1}^m\tilde k_D(z_{j-1},z_{j}):z_0,\dots,z_{m}\in D,
z_0=z,z_m=w\},$$ then
$$k_D(z,w)=\inf_m k_D^{(m)}(z,w).$$

If $l_D$ is any one of the introduced functions from above, we set
$l_D^\ast=\tanh l_D.$

Recall that $D$ is said to be balanced if $\lambda z\in D$ for any
$\lambda\in\overline{\Bbb D}$ and any $z\in D.$ Denoting by $h_D$
the Minkowski function of $D,$ i.e.,
$$h_D(z)=\inf\{t>0:z/t\in D\},\ z\in\Bbb C^n,$$ then
$$D=D_h=\{z\in\Bbb C^n:h_D(z)<1\}.$$
We point out that $D$ is  pseudoconvex if and only if  $\log h_D$
is a plurisubharmonic function. Set $\widehat h_D=h_{\widehat D},$
where $\widehat D$ is the convex hull of $D.$

Let us summarize some well-known facts about relations between
$h_D,$ $\widehat h_D,$ and the functions from above, where one of
their arguments is the origin.

\begin{proposition}\label{pr1} (cf. \cite{Jar-Pfl1}) Let $D\subset\Bbb C^n$
be a balanced domain and $a\in\Bbb D.$ Then

(i) $\widehat h_D\le c_D^\ast(0,\cdot)\le\tilde
k^\ast_D(0,\cdot)\le h_D;$

(ii) $k^\ast_D(0,a)=h_D(a)\iff h_D(a)=\widehat h_D(a).$

If, in addition, $D$ is pseudoconvex, then
$$\tilde k^\ast_D(0,\cdot)=h_D.$$
\end{proposition}

Proposition \ref{pr1} shows that at a point $a\in D$ the value of
$k_D(0,\cdot)$ is maximal if and only if $D$ is "convex in the
direction of $a$", i.e., $h_D(a)=\widehat h_D(a)$. In fact, more
is true as the following result shows.\footnote{The proofs of the
following propositions and examples will be given in Section 2.}

\begin{proposition}\label{pr2} Let $D\subset\Bbb C^n$ be a balanced domain and
$a\in D.$ Then $$(k^{(3)}_D(0,a))^\ast=h_D(a)\iff h_D(a)=\widehat
h_D(a).$$
\end{proposition}

\noindent{\bf Remark.} We do not know if the number 3 can be
replaced by 2.

\smallskip

Conversely, one may ask whether the fact that $l_D(0,a)$ ($l_D$ as
above) is "minimal" (i.e., $l_D(0,a)=h_G(a)$ for some domain $G$
containing $D$) implies also some convexity property. Here we
start with the following result.

\begin{proposition}\label{pr8} Let $D\subset\Bbb C^n$ be a bounded balanced
domain and $G\subset\Bbb C^n$ be a pseudoconvex balanced domain
with $G\supset D$. Assume that $h_D$ is continuous at some $a\in
D$ and $\overline G$ contains no nontrivial analytic discs through
$a/h_G(a),$ $h_G(a)\neq 0.$ Then
$$\tilde k^\ast_D(0,a)=h_G(a)\iff h_D(a)=h_G(a).$$
\end{proposition}

\noindent{\bf Remarks.} (a) Since the envelope of holomorphy
$\mathcal E(D)$ of a balanced domain $D$ is balanced (see
\cite{Jar-Pfl2}, Remark 3.1.2(b)), one may apply the above result
for $D$ and $\mathcal E(D)$.

(b) If $h_G$ is continuous near $a$ and $\partial G$ contains no
nontrivial discs through $a/h_G(a),$ then the maximum principle
implies that $\overline G$ contains no nontrivial analytic discs
through $a/h_G(a),$ too.

(c) In light of Proposition \ref{pr8}, it is natural to ask
whether there is a non-pseudoconvex balanced domain $D$ such that
$h_D=\tilde k_D^\ast(0,\cdot)$ on $D$. The authors do not know the
answer.

\smallskip

The following examples show that the assumptions about continuity
of $h_D$ at $a$ and discs in Proposition \ref{pr8} are essential.

\begin{example}\label{ex9} If $D=\Bbb
D^2\setminus\{(t,t):|t|\ge 1/2\},$ $d=(t,t),$ $|t|<1/2,$ then
$$h_D(d)=2|t|\hbox{\ \ but\ \ }\tilde k_D^\ast(0,d)=|t|.$$
On the other hand, $\overline{\Bbb D^2}$ contains no nontrivial
analytic discs through any point of $\partial\Bbb
D\times\partial\Bbb D.$

\end{example}

Moreover, the following example gives a balanced Reinhardt domain
$D$ such that  $G=\mathcal E(D)=\widehat D$ has nontrivial
analytic discs in its boundary for which Proposition \ref{pr8}
fails to hold. Note that in this case $h_D$ and $h_G$ are
continuous functions.

\begin{example}\label{ex10} Let $0<a<1$ and
$$D=\{z\in\Bbb D^2:|z_2|^2-a^2<2(1-a^2)|z_1|\}.$$
Then $D$ is a balanced Reinhardt domain and $\mathcal E(D)=\Bbb
D^2$ (cf. \cite{Jar-Pfl2}). On the other hand, if $c=(0,d)$,
$|d|<a,$  then $$h_D(c)=\frac{|d|}{a}>\tilde
k^\ast_D(0,c)=|d|=\widehat h_D(c).$$
\end{example}

For a balanced domain $D$ and $a\in D$ set
$$\Bbb D_{D,a}:=\{\lambda a:|\lambda|h_D(a)<1\}.$$
Assuming minimality along the whole slice $\Bbb C a\cap D=\Bbb
D_{D,a}$ for some $a\in D$ we have

\begin{proposition}\label{pr3} Let $D\subset\Bbb C^n$ be a taut balanced
domain, $a\in D$ and $m\in\Bbb N.$ Then
$$(k^{(m)}_D(0,\tilde a))^\ast=\widehat h_D(\tilde a),\ \tilde a\in\Bbb
D_{D,a} \iff h_D(a)=\widehat h_D(a).$$
\end{proposition}

Recall that a domain $D\subset\Bbb C^n$ is said to be taut if
$\mathcal O(\Bbb D,D)$ is a normal family. Note that a balanced
domain $D$ is taut if and only if $h_D$ is a continuous
plurisubharmonic function and $(h_D)^{-1}(0)=0.$ Weakening the
continuity assumption for $h_D$ the following statement remains
true.

\begin{proposition}\label{pr4} Let $D\subset\Bbb C^n$ be a bounded
pseudoconvex balanced domain, $a\in D$ and $m\in\Bbb N.$ Assume
that $h_D$ is continuous at $a$ and $\partial D$ contains no
nontrivial analytic discs through $a/h_D(a).$ Then
$$(k^{(m)}_D(0,\tilde a))^\ast=\widehat h_D(\tilde a),\quad \tilde a\in\Bbb D_{D,a}
\iff h_D(a)=\widehat h_D(a).$$
\end{proposition}

\noindent{\bf Remark.} We do not know if the condition about discs
is superfluous or not. On the other hand, continuity and
pseudoconvexity are essential as Examples \ref{ex9} and \ref{ex10}
have shown (see also Example \ref{ex6} below).
\smallskip

Recall that a boundary point $b$ of a domain $D$ in $\Bbb C^n$ is
said to be a \textit{local weak barrier point} if there are a
neighborhood $U$ of $b$ and a negative plurisubharmonic function
$u$ on $D\cap U$ such that $\lim_{D\ni z\to b}u(z)=0.$

Proposition \ref{pr4} is a consequence of the following

\begin{proposition}\label{pr5} Let $b$ be a local weak barrier point of a
bounded domain $D\subset\Bbb C^n.$ If $\partial D$ contains no
nontrivial analytic discs through $b,$ then
$$\lim_{w\to b}k^{(m)}_D(z,w)=\infty,\ z\in D,\ m\in\Bbb N.$$
\end{proposition}

We point out that the assumption about tautness is essential in
Proposition \ref{pr3} as the following examples show.

\begin{example}\label{ex6} The unit ball $\Bbb B$ in $\Bbb C^2$
contains a proper non-taut pseudoconvex balanced domain $D$ with
$\widehat D=\Bbb B$ such that
$$(k^{(2)}_D(0,\cdot))^\ast=||\cdot||.$$
\end{example}

\begin{example}\label{ex7} There is an unbounded pseudoconvex balanced
Reinhardt domain $D$ in $\Bbb C^2$ with continuous Minkowski
function and a point $a\in D$ such that $\widehat h_D(a)>0,$
$\partial D$ contains a nontrivial analytic disc through
$a/h_D(a),$
$$(k^{(2)}_D(0,\lambda a))^\ast=|\lambda|\widehat h_D(a),\ |\lambda|\le 1,$$
but even
$$c^\ast_D(0,\lambda a)>|\lambda|\widehat h_D(a),\ 1<|\lambda|<1/h(a).$$
\end{example}

\noindent{\bf Remark.} Despite of these examples, we do not know
any example of a taut balanced domain $D$ such that $(k^{(m)}_D(0,
a))^\ast=\widehat h_D(a)$ for some $m\in\Bbb N$ and some $a\in D,$
but $h_D(a)>\widehat h_D(a).$

\section{Proofs}

\begin{proof}[Proof of Proposition \ref{pr2}] We have only
to prove that $$(k^{(3)}_D(0,a))^\ast=h_D(a)\Rightarrow
h_D(a)=\widehat h_D(a).$$

First, we shall show that
\begin{equation}\label{1}(k^{(2)}_D(0,\lambda a))^\ast=|\lambda|h_D(a),\
\lambda\in\Bbb D.
\end{equation}
We may assume that $h_D(a)\neq 0.$ Taking the disc $\Bbb D\ni t\to
ta/h_D(a)$ as a competitor for $\tilde k_D(\lambda a,a)$ gives
$$\tilde k_D(\lambda a,a)\le p(h_D(\lambda a),h_D(a)),$$ where $p$
denotes the Poincar\'e distance. This and the inequality
$$p(0,h_D(a))=k^{(3)}_D(0,a)\le k^{(2)}_D(0,\lambda a)+\tilde k_D(\lambda
a,a)$$ imply that
$$p(0,|\lambda|h_D(a))=p(0,h_D(a))-p(|\lambda|h_D(a),h_D(a))\le
k^{(2)}_D(0,\lambda a).$$ So $$(k^{(2)}_D(0,\lambda
a))^\ast\ge|\lambda|h_D(a)$$ and the opposite inequality always
holds.

It follows from (\ref{1}) that
$$\lim_{\lambda\to 0}\frac{k^{(2)}_D(0,\lambda
a)}{|\lambda|}=h_D(a).$$ On the other hand, by Proposition 2 in
\cite{Nik-Pfl2}, this limit does not exceed
$$\kappa^{(2)}_D(0;a):=\inf\{\kappa_D(0;a_1)+\kappa_D(0;a_2):a_1+a_2=a\},$$ where
$\kappa_D$ denotes the Kobayashi--Royden pseudometric of $D.$
Since $\kappa_D(0;\cdot)\le h_D$ (cf. \cite{Jar-Pfl1}), we
conclude that
$$h_D(a)\le h_D(a_1)+h_D(a_2)\hbox{ if }a_1+a_2=a,$$
which means that $h_D(a)=\widehat h_D(a).$
\end{proof}

\begin{proof}[Proof of Proposition \ref{pr8}]
We have only to prove that
$$\tilde k^\ast_D(0,a)=h_G(a)\Rightarrow h_D(a)\le h_G(a).$$
Let $(\varphi_j)\subset\mathcal O(\Bbb D,D)$ and $\alpha_j\to
h_G(a)$ be such that $\varphi_j(0)=0$ and $\varphi_j(\alpha_j)=a.$
Writing $\varphi_j$ in the the form
$\varphi_j(\lambda)=\lambda\psi_j(\lambda),$ $\psi_j\in\mathcal
O(\Bbb D,\Bbb C^n),$ it follows by the maximum principle that
$h_G\circ\psi_j\le 1$ and hence $\psi_j\in\mathcal O(\Bbb
D,\overline G).$ Since $D$ is bounded, we may assume that
$\varphi_j\to\varphi\in\mathcal O(\Bbb D,\overline D)$ and then
$\psi_j\to\psi\in\mathcal O(\Bbb D,\overline G).$ On the other
hand, since $\overline G$ contains no nontrivial analytic discs
through $\psi(h_G(a))=b,$ where $b=a/h_G(a),$ it follows that
$\psi(\Bbb D)=b.$ Using that $h_D$ is continuous at $b,$ we get
that
$$1>h_D(\varphi_j(\lambda))\to|\lambda|h_D(b),\ \lambda\in\Bbb D.$$
Letting $\lambda\to 1$ leads to $h_D(b)\le 1$ which is the desired
inequality.
\end{proof}

\begin{proof}[Proof of Example \ref{ex9}] We have only to prove
that
$$\tilde k^\ast_D(0,d)\le|t|.$$
For any $r\in(|t|,1)$ we may choose $\alpha\in\Bbb D$ such that
$t=\varphi(t/r),$ where
$\varphi(\lambda)=\lambda\frac{\lambda-\alpha}{1-\overline{\alpha}\lambda}.$
Then $\psi=(r\mbox{id},\varphi)\in\mathcal O(\Bbb D,D)$ is a
competitor for $\tilde k^\ast_D(0,d)$ which shows that $\tilde
k^\ast_D(0,d)\le|t|/r.$ It remains to let $r\to 1.$
\end{proof}

\noindent{\bf Remarks.} (a) Let $D$ be the domain from Example
\ref{ex9}. Note that even
$$\tilde k_D(0,\cdot)=\tilde k_{\Bbb D^2}(0,\cdot).$$ It is enough
to prove that $\tilde k_D(0,a)\le|a_1|$ for $a=(a_1,a_2)\in D,$
$a_1\neq a_2,$ $|a_1|\ge|a_2|.$ For this, take the discs
$\psi(\lambda)=(\lambda,\lambda a_2/a_1)$ as a competitor for
$\tilde k_D(0,a).$

On the other hand, if $a_1=(0,b)$ and $a_2=(b,0),$ $b\in\Bbb D,$
then
$$\tilde k_D(a_1,a_2)=\tilde k_{\Bbb D^2}(a_1,a_2)\iff |b|\le 4/5.$$

Indeed, using the M\"obius transformation
$\psi_b(\lambda)=\frac{\lambda-b}{1-\overline{b}\lambda},$ we get
that
$$\tilde k_D(a_1,a_2)=\tilde k_{D_b}(0,a),$$
where $a=(b,-b)$ and $D_b=\Bbb
D^2\setminus\{(\psi_b(\lambda),\lambda):1/2\le |\lambda|<1\}.$

For $|b|<4/5,$ it is easy to check that
$\varphi=(\mbox{id},-\mbox{id})\in\mathcal O(\Bbb D, D_b).$ This
implies that $\tilde k^\ast_{D_b}(0,a)\le|b|$ and hence $\tilde
k_D(a_1,a_2)=\tilde k_{\Bbb D^2}(a_1,a_2).$

To get the same for $|b|=4/5,$ it is enough to take $r\varphi,$
$r\in(0,1),$ as a competitor for $\tilde k^\ast_{D_b}(0,a)$ and
then to let $r\to 1.$

Assume now that $|b|>4/5$ and $\tilde k_D(a_1,a_2)=\tilde k_{\Bbb
D^2}(a_1,a_2).$ Then we may find discs $\varphi_j\in\mathcal
O(\Bbb D,D_b)$ such that $\varphi_j(0)=0$ and
$\varphi_j(\alpha_j)=a,$ where $\alpha_j\to b.$ It follows by the
Schwarz-Pick lemma that $\varphi_j\to\varphi.$ On the other hand,
$\varphi(\Bbb D)\cap\{(\psi_b(\lambda),\lambda):1/2<|\lambda|<1\}$
is a singleton which contradicts to Hurwitz's theorem.

(b) Note that (a) shows that Theorem 3.4.2 in \cite{Jar-Pfl1} is
in some sense sharp. On the other hand, this theorem implies that,
if $D_n=\Bbb D^n\setminus\{(t,\dots,t):|t|\ge 1/2\},$ $n\ge 3,$
then
$$\tilde k_{D_n}=\tilde k_{\Bbb D^n}.$$

\begin{proof}[Proof of Example \ref{ex10}] We have only to prove
that $\tilde k^\ast_D(0,c)\le d,$ $d\in(0,a).$  For this, it is
enough to show that $\varphi=(\psi,\mbox{id})\in\mathcal O(\Bbb
D,D),$ where $\psi(\lambda)=\lambda\frac{\lambda-d}{1-d\lambda}.$
Since $|\psi(\lambda)|\ge x\frac{x-d}{1-xd}$ for $x=|\lambda|,$ we
have to check that $$(x^2-a^2)(1-dx)<2(1-a^2)x(x-d),\hbox{
i.e.,}$$
$$dx^3+(1-2a^2)x^2-d(2-a^2)x+a^2>0.$$
Note that this inequality is true for $x=0$. Using that $x\in
(0,1)$ and $d\in(0,a)$, it suffices to prove that
$$ax^3+(1-2a^2)x^2-a(2-a^2)x+a^2 \ge 0$$
which is equivalent to the obvious inequality $(x-a)^2(ax+1)\ge
0.$
\end{proof}

\begin{proof}[Proof of Proposition \ref{pr3}]
It is enough to show that if $\Bbb D_{D,a}\ni a_k\to a/h_D(a),$
$a\neq 0,$ and $(k^{(m)}_D(0,\tilde a_k))^\ast=\widehat h_D(\tilde
a_k),$ then $h_D(a)=\widehat h_D(a).$ For this, recall that if $D$
is a taut domain, then (cf. \cite{Jar-Pfl1}, Proposition 3.2.1)
$$\lim_{w\to\partial D}k^{(m)}_D(z,w)=\infty,\ w\in D.$$ In our case
this implies that $\widehat h_D(\tilde a_k)\to 1$ and hence
$\widehat h_D(a/h(a))=1.$
\end{proof}

\begin{proof} [Proof of Proposition \ref{pr5}]

We shall proceed by induction on $m.$ For this, we shall need the
following.

\begin{lemma}\label{l11} Under the assumptions of Proposition \ref{pr5},
for any $(\varphi_k)\subset\mathcal O(\Bbb D,D)$ with $\varphi_k
(0)\to b$ one has that $\varphi_k\to b$ locally uniformly on $\Bbb
D.$
\end{lemma}

Assuming Lemma \ref{l11} easily implies that Proposition \ref{pr5}
is true for $m=1.$ Suppose that this statement is true for some
$m-1\in\Bbb N$ but false for $m.$ Then we may find $z\in D$ and
sequences $(z_{j,k})_k\subset D,$ $0\le j\le m,$ such that
$z_{0,k}=z,$ $z_{m,k}\to b$ and
$$\sup_k\sum_{j=1}^m\tilde k_D(z_{j-1,k},z_{j,k})<\infty.$$
In virtue of our induction hypothesis one has that
$z_{m-1,k}\not\to b.$ Passing to a subsequence, we may assume that
$z_{m-1,k}\to a\in\overline D,$ $a\neq b$ and $(\tilde
k_D(z_{m-1,k},z_{m,k}))^\ast \to r<1.$ Then there are
$\varphi_k\in\mathcal O(\Bbb D,D)$ with $\varphi_k(0)=z_{m,k},$
$\varphi_k(r_k)=z_{m-1,k}$ and $r_k\to r$ which contradicts Lemma
11.
\end{proof}

\begin{proof} [Proof of Lemma \ref{l11}] Since $D$ is bounded, it
is easily seen that for any neighborhood $U$ of b there is another
neighborhood $V\subset U$ of $b$ and a number $s\in(0,1]$ such
that, if $\varphi\in\mathcal O(\Bbb D,D)$ with $\varphi(0)\in V,$
then $\varphi(s\Bbb D)\subset D\cap U.$

Now we choose $U$ such that there is a negative plurisubharmonic
function $u$ on $D\cap U$ with $\lim_{D\ni z\to b} u(z)=0.$
Applying Proposition 3.4 in \cite{Nik-Pfl1} implies that $b$ is a
{\it $t$-point} for $D\cap D$ which means that

{\it $(\ast)$ $(\varphi_j|_{s\Bbb D})$ is compactly divergent
w.r.t. $D\cap U$.}

Assume now that $(\varphi_j)$ does not converge to $b$. Passing to
a subsequence, we may assume that $\varphi_j\to\varphi\in\mathcal
O(\Bbb D,\overline{D})$ and $\varphi(\Bbb D)\neq\{b\}.$ It follows
by $(\ast)$ that $\varphi(s\Bbb D)\subset\partial(D\cap U).$ Since
$\varphi(0)=b\in\partial D\cap U,$ we may find $s'\in(0,s]$ such
that $\varphi(s'\Bbb D)\subset\partial D.$ Using that $\partial D$
contains no nontrivial analytic discs through $b,$ we get that
$\varphi(s'\Bbb D)=\{b\}$ which contradicts the identity
principle.

\end{proof}

\begin{proof}[Proof of Example \ref{ex6}]
J. Siciak (cf. \cite{Jar-Pfl1}, Example 3.1.12) constructed a
plurisubharmonic function $\psi:\Bbb C^2\to[0,\infty)$ such that
$\psi(\lambda z)=|\lambda|\psi(z)$ ($\lambda\in\Bbb C,z\in\Bbb
C^n$), $\psi\not\equiv 0$, but $\psi=0$ on a dense subset $S$ of
$\Bbb C^n.$ Set
$$D=\{z\in\Bbb C^2:||z||+\psi(z)<1\}.$$
For any $a\in D$ we may choose a sequence $S\cap D\supset(z_k)\to
a.$ Then
$$\tanh^{-1}||a||=k^{(2)}_{\Bbb B_2}(0,a)\le k^{(2)}_D(0,a)\le
\tilde k_D(0,z_k)+\tilde k_D(z_k,a).$$ Letting $k\to\infty$ gives
$k^{(2)}_D(0,a)=\tanh^{-1}||a||.$
\end{proof}

\begin{proof}[Proof of Example \ref{ex7}] Setting
$$D=\{z\in\Bbb C^2:|z_1|<1,|z_1^{34}z_2^{55}|<1\},$$
then
$$c_D^\ast(0,z)=\max\{|z_1|,|z_1z_2|,|z_1^2z_2^3|,|z_1^5z_2^8|,|z_1^{13}z_2^{21}|,
|z_1^{34}z_2^{55}|\}$$ (see \cite{Jar-Pfl1}, Example 2.7.12).
Since $\widehat D=\Bbb D\times\Bbb C,$ we have that
$$c_D^\ast(0,z)>\widehat h_D(z)=|z_1|\iff |z_2|>1.$$
On the other hand, note that $\{0\}\times\Bbb C\subset D$ and
$\Bbb D\times\{z_2\}\subset D,$ $|z_2|\le 1.$ So if
$|z_1|<1,|z_2|\le 1,$ then
$$k^{(2)}_D(0,(z_1,z_2))\le k_D(0,(0,z_2))+k_D((0,z_2),(z_1,z_2))\le
k_{\Bbb D}(0,z_1)$$ and hence
$$(k^{(2)}_D(0,(z_1,z_2)))^\ast\le|z_1|=\widehat h_D((z_1,z_2)).$$ Since the
opposite inequality always holds and $\partial D$ contains a
nontrivial analytic disc through any point $b=(b_1,b_2)\in\partial
D$ with $|b_1|\neq|b_2|,$ it follows that any point $a=(a_1,a_2)$
with $|a_1|<|a_2|=1$ has the desired properties.
\end{proof}

\end{document}